\documentclass{article}

\usepackage[utf8]{inputenc} 
\usepackage[T1]{fontenc}    
\usepackage{hyperref}       
\usepackage{url}            
\usepackage{booktabs}       
\usepackage{amsfonts}       
\usepackage{nicefrac}       
\usepackage{microtype}      
\usepackage{lipsum}
\usepackage{fancyhdr}       
\usepackage{graphicx}       
\graphicspath{{media/}}     
\usepackage{amsmath}
\usepackage{longtable}
\usepackage{makecell}
\usepackage{xcolor}
\usepackage{float}

\title{Partially invariant solution with an arbitrary surface of blow-up for the gas dynamics equations admitting pressure translation
}

\author{
	Dilara Siraeva \\
	Department of Mathematics\\ North Carolina State University, Raleigh\\
	\texttt{sirdilara@gmail.com} 
}

\begin{document}
\maketitle
\begin{abstract}
We applied a method of symmetry reduction to the gas dynamics equations with a special form of the equation of state. This equation of state is a pressure represented as the sum of a density and an entropy functions. The symmetry Lie algebra of the system is 12--dimensional. One, two and three--dimensional subalgebras were considered. In this article, four--dimensional subalgebras are considered for the first time. Specifically, invariants are calculated for 50 four--dimensional subalgebras. Using invariants of one of the subalgebras, a symmetry reduction of the original system is calculated. The reduced system is a partially invariant submodel because one gas-dynamic function cannot be expressed in terms of the invariants. The submodel leads to two families of exact solutions, one of which describes the isochoric motion of the media, and the other solution specifies an arbitrary blow-up surface. For the first family of solutions, the particle trajectories are parabolas or rays; for the second family of solutions, the particles move along cubic parabolas or straight lines. From each point of the blow-up surface, particles fly out at different speeds and end up on a straight line at any other fixed moment in time. A description of the motion of particles for each family of solutions is given.
\end{abstract}

\noindent {\bf Keywords}: Gas dynamics equations, state equation, submodel, invariant, exact solution, blow-up.

\noindent {\bf MSC 2020:} 35B06, 35C99.

\section{Introduction}
The gas dynamics equations are~\cite{Ovs1}
\begin{equation}\label{siraeva:eq0}
	\begin{array}{c}
		D\vec{u}+\rho^{-1}{\bf\nabla} P=0,\\[2mm] D\rho+\rho\,{\rm div}\vec{u}=0,\\[2mm] DP+\rho f_{\rho}{\rm div}\vec{u}=0,
	\end{array}
\end{equation}
where $D=\partial_t+(\vec{u}\cdot{\bf\nabla})$ is a total differentiation operator; $t$ is time; $\nabla=\partial_{\vec{x}}$ is a gradient by spatial independent variables $\vec{x}$;   $\vec{u}$ is a velocity vector; $\rho$ is a density; $P$ is a pressure. An arbitrary equation of state has the form
\begin{equation}\label{sir:arbitrary_eq_state}
	\begin{array}{c}
P=f(\rho,S).
\end{array}
\end{equation}
An entropy $S$  is determined from the equality \eqref{sir:arbitrary_eq_state}. The last equation of system \eqref{siraeva:eq0} can be replaced by an equation for the entropy
\[DS=0.\]
Finding new exact solutions to equations~\eqref{siraeva:eq0}, \eqref{sir:arbitrary_eq_state} is a non-trivial problem. New solutions make it possible to describe the motion of particles, to test numerical calculations of gas-dynamic problems, and to create new numerical schemes for calculating special boundary value problems. The symmetry analysis of differential equations~\cite{Sir:lit:Ovs2,Sir:Olver1986,Sir:Bluman1989, Sir:Ibrag1994,Sir:Ibrag2001,Sir:Hydon2000,Sir:Arrigo2015} allows us to obtain new exact solutions to the gas dynamics equations. The application of this method implies obtaining a reduction of the original system, which is easier to solve. The method is based on the use of a transformation group with respect to which the gas dynamics equations remain invariant. The transformation group of system~\eqref{siraeva:eq0},~\eqref{sir:arbitrary_eq_state} involves time translation, space translations, Galilean translations, rotations and uniform dilatation of space-time. In the 20th century, a Soviet mathematician Lev Ovsyannikov calculated that for some special equations of state, this group has additional transformations~\cite{Ovs1}. In the case of the special equation of state  \cite{Ovs1}
\begin{equation}\label{Siraeva:US_special}
	\begin{array}{c}
		P=f(\rho)+S
	\end{array}
\end{equation}
the transformation group allows for the additional transformation of pressure translation. Such transformation group corresponds to a 12--dimensional Lie algebra $L_{12}$  whose basis is given in the next section.

For each subalgebra of $L_{12}$, one can calculate invariants --- functions that are unchanged by the group actions and annihilated by the infinitesimal generators. The following terminology will be used in this article. Invariant submodel is system~\eqref{siraeva:eq0},~\eqref{Siraeva:US_special} which is written in terms of invariants of subalgebra. Rank is a number of independent variables of the submodel. Submodel is partially invariant if one or more gas-dynamic functions are not expressed through invariants. The number of such functions is a defect of the submodel. 

All subalgebras of the Lie algebra $L_{12}$, up to internal automorphisms, were listed in the optimal system of nonsimilar subalgebras in~\cite{siraeva2014optimal}. Let us list some previously obtained results for two and three--dimensional subalgebras. Two 2--dimensional subalgebras of Lie algebra $L_{12}$ define partially invariant submodels of rank 3 and defect 1. Their reductions to invariant submodels were obtained in~\cite{SirMnSys}. Invariant submodels of rank 2 in the canonical form of Lie algebra $L_{12}$ were constructed in~\cite{Sir16,SirSibJIM1}. An example of a description of particles motion for solution from invariant submodel of rank 2 is given in~\cite{Sir4Khab5}. Exact solutions were obtained in \cite{Sir:lit:Sir1JoPh20, Sir:lit:SEMR2021,Sir:lit:MDPI2022} for four submodels of rank 1.

The present article is dedicated to four--dimensional subalgebras of Lie algebra $L_{12}$, which were not considered. Four--dimensional subalgebras were considered for other Lie algebras. For instance, for all 48 types of four--dimensional subalgebras of Lie algebra $L_{11}$ the simplest partially invariant solutions of rank 1 and defect 1 were considered~\cite{Sir:lit:Khabirov2019}. Eight simple invariant solutions of rank 0 were obtained from four--dimensional subalgebras containing the projective generator in the case of a system of gas dynamics equations with the state equation of a monatomic gas. The motion of gas particles as a whole was constructed for the isentropic solution. The solutions obtained have a density singularity on a constant or moving plane, which is a boundary with vacuum or a wall~\cite{Sir:lit:Nikon2023}.
	
\section{The system and its symmetry Lie algebra}
The present article considers system~\eqref{siraeva:eq0},~\eqref{Siraeva:US_special} in the Cartesian coordinate system. We have
\begin{equation*}
	\begin{array}{c}
	\vec{x}=x\vec{i}+y\vec{j}+z\vec{k},\\[2mm] \nabla=\vec{i}\partial_x+\vec{j}\partial_y+\vec{k}\partial_z,\\[2mm] \vec{u}=u\vec{i}+v\vec{j}+w\vec{k},
	\end{array}
\end{equation*}
where $\vec{i}$, $\vec{j}$, $\vec{k}$ is an orthonormal basis.

It is known that the system \eqref{siraeva:eq0} with an arbitrary equation of state~\eqref{sir:arbitrary_eq_state} admits an 11--dimensional Lie algebra  $L_{11}$. The equations of state for which an extension of the 11--dimensional Lie algebra occurs were listed in \cite{Ovs1}. All non-isomorphic Lie algebras of the group classification according to the equation of state were given in~\cite{Khab1}. For each Lie algebra, the method of calculating nonsimilar subalgebras was finally completed in~\cite{Khab2}. Invariant submodels of rank 3, 1 \cite{Chir1Khab4,Sir:lit:Khab3}, and 2 \cite{Sir:lit:Mamont1,Sir:lit:Mamont2}  were constructed for an 11--dimensional Lie algebra.

The equations~\eqref{siraeva:eq0} with an arbitrary equation of state \eqref{sir:arbitrary_eq_state} are invariant under the action of the Galilean group extended by uniform dilatation~\cite{Chir1Khab4}:
\begin{equation}\label{Sir:GroupGal_rashir}
	\begin{array}{l}
		1^o\!.\; \vec{x}^{\,*}=\vec{x}+\vec{a} \text{ (space translations)},\\
		
		2^o\!.\; t^{*}=t+a_0 \text{ (time translation)},\\
		
		3^o\!.\; \vec{x}^{\,*}=O\vec{x}, \vec{u}^{\,*}=O\vec{u}, OO^{T}=E, \det O=1 \text{ (rotations)},\\
		
		4^o\!.\; \vec{x}^{\,*}=\vec{x}+t\vec{b}, \vec{u}^{\,*}=\vec{u}+\vec{b} \text{ (Galilean translations)},\\
		
		5^o\!.\; t^{*}=ct, \vec{x}^{\,*}=c\vec{x} \text{ (uniform dilatation)}.
	\end{array}
\end{equation}
The equations~~\eqref{siraeva:eq0} with special equation of state \eqref{Siraeva:US_special}  are invariant under the action of transformations \eqref{Sir:GroupGal_rashir} and the pressure translation~\cite{Chir1Khab4}:
\begin{equation}\label{Sir:peren_po_p}
	\begin{array}{l}
		6^o\!.\;P^{*}=P+P_0.
	\end{array}
\end{equation}

Each of the transformations~\eqref{Sir:GroupGal_rashir},~\eqref{Sir:peren_po_p} maps a solution of the system~\eqref{siraeva:eq0}, \eqref{Siraeva:US_special} to another solution. Therefore, in what follows, the solutions are considered up to the transformations \eqref{Sir:GroupGal_rashir}, \eqref{Sir:peren_po_p}. 

The 12--dimensional Lie algebra $L_{12}$ corresponds to the transformation group \eqref{Sir:GroupGal_rashir}, \eqref{Sir:peren_po_p}. The basis generators of $L_{12}$ in the Cartesian coordinate system have the form~\cite{Ovs1}
\begin{equation}\label{sir:infinit_generators}
	\begin{array}{c}
		X_1=\partial_{x}, \quad X_2=\partial_{y}, \quad X_3=\partial_{z},\\[2mm]  X_4=t\partial_{x}+\partial_{u}, \quad X_5=t\partial_{y}+\partial_{v}, \quad X_6=t\partial_{z}+\partial_{w},\\[2mm]  X_7=y\partial_{z}-z\partial_{y}+v\partial_{w}-w\partial_{v},\quad 
		X_8=z\partial_{x}-x\partial_{z}+w\partial_{u}-u\partial_{w},\\[2mm] X_9=x\partial_{y}-y\partial_{x}+u\partial_{v}-v\partial_{u}, \quad X_{10}=\partial_{t},\\[2mm]
		X_{11}=t\partial_{t}+x\partial_{x}+y\partial_{y}+z\partial_{z}, \quad Y_{1}=\partial _P.
	\end{array}
\end{equation}
The commutators of infinitesimal generators \eqref{sir:infinit_generators} are given in 
Table~\ref{Diss:tabl_komm:1} \cite{siraeva2014optimal}, where instead of generators $X_i$, $i=1...11$, we simply write indices $i$.
\begin{table}[h!]
	\caption{\mbox{The table of commutators of infinitesimal generators of Lie algebra~$L_{12}$}.}\label{Diss:tabl_komm:1}
	\begin{center}
		\begin{tabular}{ c p{2.4cm} p{2.4cm} p{2.4cm} p{1.7cm} p{0.2cm} }
			\hline
			{} & $\hspace{0.3cm}{1}$ $\hspace{0.5cm}{2}$ $\hspace{0.5cm}{3}$     & ${\hspace{0.3cm}4}$ ${\hspace{0.5cm}5}$ ${\hspace{0.5cm}6}$   &  $\hspace{0.3cm}{7}$  $\hspace{0.5cm}{8}$ $\hspace{0.5cm}{9}$ & $\hspace{0.2cm}{10}$ $\hspace{0.4cm}{11}$ & $Y_1$\\   \hline
			
			{1} &                 &                 & $\hspace{0.85cm}{-3}$  $\hspace{0.4cm}{2}$     &  $\hspace{1.3cm}{1}$  &\\
			
			{2} &                 &                 & $\hspace{0.3cm}{3}$    $\hspace{1cm}{-1}$    &  $\hspace{1.3cm}{2}$ & \\
			
			{3} &                 &                 & $\hspace{0.03cm}{-2}$   $\hspace{0.4cm}{1}$     &  $\hspace{1.3cm}{3}$  &\\  
			
			{4} &                 &                 & $\hspace{0.85cm}{-6}$  $\hspace{0.4cm}{5}$     &  $\hspace{0.09cm}{-1}$ & \\
			
			{5} &                 &                 & $\hspace{0.3cm}{6}$    $\hspace{0.95cm}{-4}$    &  $\hspace{0.09cm}{-2}$  &\\
			
			{6} &                 &                 & $\hspace{0.05cm}{-5}$   $\hspace{0.4cm}{4}$      & $\hspace{0.09cm}{-3}$               &  \\   
			
			{7} &  $\hspace{0.85cm}{-3}$ $\hspace{0.4cm}{2}$ & $\hspace{0.78cm}{-6}$ $\hspace{0.47cm}{5}$ & $\hspace{0.85cm}{-9}$ $\hspace{0.4cm}{8}$   &   & \\
			
			{8} &  $\hspace{0.3cm}{3}$ $\hspace{0.95cm}{-1}$  & $\hspace{0.3cm}{6}$ $\hspace{0.95cm}{-4}$ &  $\hspace{0.3cm}{9}$ $\hspace{0.95cm}{-7}$   &   & \\
			
			{9} &  $\hspace{0.03cm}{-2}$ $\hspace{0.5cm}{1}$  &    $\hspace{0.03cm}{-5}$ $\hspace{0.45cm}{4}$ & $\hspace{0.045cm}{-8}$ $\hspace{0.5cm}{7}$ & &\\ 
			
			{10}&                  &   $\hspace{0.3cm}{1}$  $\hspace{0.47cm}{2}$ $\hspace{0.47cm}{3}$ &           & $\hspace{1.1cm}{10}$ & \\
			{11}& $\hspace{0.03cm}{-1}$ $\hspace{0.09cm}{-2}$ $\hspace{0.09cm}{-3}$ &  &  &  $\hspace{0cm}{-10}$ & \\   \hline
			$Y_1$& &  &  &   & \\ \hline
		\end{tabular}
	\end{center}
\end{table}

Lie algebra $L_{12}$ is a direct sum of two ideals $L_{12}=L_{11} \oplus \{Y_1\}$. An outer automorphism of Lie algebra $L_{12}$ is $y^{o\,*}=b_{1}\,y^{o}$ and internal automorphisms are in Table~\ref{siraeva:table_autom}, where $x^i$, $x^{i*}$, $y^o$, $y^{o*}$ are parameters of subalgebras $X=x^iX_i+y^oY_1$, $X^{*}=x^{i*}X_i+y^{o*}Y_1$; notations for projections are $p_1(x)=(x^1,x^2,x^3)$, $p_2(x)=(x^4,x^5,x^6)$, $p_3(x)=(x^7,x^8,x^9)$, $p_1(x^{*})=(x^{1*},x^{2*},x^{3*})$, $p_2(x^{*})=(x^{4*},x^{5*},x^{6*})$; while $\vec{\alpha}_{1}=(a_{1},a_{2}, a_{3})$, $\vec{\alpha}_{2}=(a_{4}, a_{5}, a_{6})$, $a_{10}$, $a_{11}$, $b_{1}$ are the parameters of automorphisms, $O$ is the rotation matrix defined by the rotation angles around one of the orthogonal axes, and $\varepsilon_{1}, \varepsilon_{2}$ are noticed discrete automorphisms.
\begin{table}[h!]
  \begin{center}
            \caption{Internal automorphisms of Lie algebra $L_{12}$}\label{siraeva:table_autom}
		\begin{tabular}{|c|p{9.5cm}|}
			\hline
			T & $p_{1}(x^{*})=p_{1}(x)+x^{11}\vec{\alpha_{1}}-\vec{\alpha_{1}}\times p_{\,3}(x)$ \\ \hline
			 $\Gamma$ & $p_{1}(x^{*})=p_{1}(x)-x^{10}\vec{\alpha_{2}},\ \ p_{\,2}(x^{*})=p_{\,2}(x)-\vec{\alpha_{2}}\times p_{\,3}(x)$ \\ \hline
			O & $p_{1}(x^{*})=Op_{1}(x),\ \ p_{\,2}(x^{*})=Op_{\,2}(x),\ \ p_{\,3}(x^{*})=Op_{\,3}(x)$ \\ \hline
			$A_{10}$ &  $p_{1}(x^{*})=p_{1}(x)+a_{10}p_{\,2}(x),\ \ x^{10*}=x^{10}+a_{10}x^{11}$ \\ \hline
			$A_{11}$ & $p_{1}(x^{*})=a_{11}p_{1}(x),\ \ x^{10*}=a_{11}x^{10}$\\ \hline
			$\varepsilon_{1}$ & $p_{1}(x^{*})=-p_{1}(x),\ \ p_{\,2}(x^{*})=-p_{\,2}(x)$ \\ \hline
			$\varepsilon_{2}$ & $p_{\,2}(x^{*})=-p_{\,2}(x),\ \ x^{10*}=-x^{10}$ \\ \hline
		\end{tabular} 
  \end{center}
		\end{table}

\section{Invariants of four--dimensional subalgebras of Lie algebra $L_{12}$}\label{siraeva:section3}
Invariants of 50 four--dimensional subalgebras of Lie algebra $L_{12}$ from the optimal system of subalgebras~\cite{siraeva2014optimal} are listed below. When calculating the invariants, the software "Maple" was used. The bases of the subalgebras are represented in the Cartesian coordinate system. The density $\rho$ is an invariant for any subalgebra of Lie algebra $L_{12}$. The number of the subalgebra from the optimal system of subalgebras~\cite{siraeva2014optimal} with parameters is indicated on the left.
\begin{enumerate}  
  
  \item[ $\begin{array}{ll}
    4.18
  \end{array}
  $] 
  
  $\begin{array}{cc}
    \text{Basis}\colon & X_2,\ \ X_3,\ \ X_7,\ \ Y_{1}+aX_4+X_{11}\\ 
    \text{Inv.:}&\dfrac{x}{t}-a\ln|t|,\ \ u-\dfrac{x}{t},\ \ v^2+w^2,\ \ P-\ln|t|\\
  \end{array}
  $
  
  \item[$\begin{array}{cc}
    4.25\\
    e\neq 0
  \end{array}
  $] 
     $\begin{array}{ll}
    \text{Basis}\colon & X_1,\ \ aX_4+X_5,\ \ bX_4+cX_6+X_{11},\ \ Y_{1}+dX_4+eX_6,\ \ d^2+e^2=1\\ 
    \text{Inv.:}& u-a\dfrac{y}{t}-\dfrac{dz}{et}+\left(\dfrac{cd}{e}-b\right)\ln|t|,\ \ v-\dfrac{y}{t},\ \ w-\dfrac{z}{t},\ \ P-\dfrac{z}{et}+\dfrac{c}{e}\ln|t|
    \end{array}
  $
  
  \item[$\begin{array}{ll}
  4.25\\
  e=0,\\d=1
  \end{array}$]
 
    $\begin{array}{ll}
    \text{Basis}\colon & X_1,\ \ aX_4+X_5,\ \ bX_4+cX_6+X_{11},\ \ Y_{1}+X_4\\ 
    \text{Inv.:}& \dfrac{z}{t}-c\ln|t|,\ \ v-\dfrac{y}{t},\ \ w-\dfrac{z}{t},\ \ P-u+a\dfrac{y}{t}+b\ln|t|
    \end{array}
  $
  
  \item[$\begin{array}{ll}
  4.26\\
  b\neq 0,\\ c\neq 0
  \end{array}$]
 
    $\begin{array}{ll}
    \text{Basis}\colon & X_1,\ \ X_4,\ \ aX_5+X_{11},\ \ Y_{1}+bX_5+cX_6,\ \ b^2+c^2=1,\ \ a\neq0\\ 
    \text{Inv.:}& \dfrac{z}{ct}-\dfrac{y}{bt}+\dfrac{a}{b}\ln|t|,\ \ v-\dfrac{y}{t},\ \ w-\dfrac{z}{t},\ \ P-\dfrac{z}{ct}
    \end{array}
  $
  
   \item[$\begin{array}{ll}
  4.26\\
  b=0\\ c=1
  \end{array}$]
 
    $\begin{array}{ll}
    \text{Basis}\colon & X_1,\ \ X_4,\ \ aX_5+X_{11},\ \ Y_{1}+X_6,\ \ a\neq0  \\ 
    \text{Inv.:}& \dfrac{y}{t}-a\ln|t|,\ \ v-\dfrac{y}{t},\ \ w-\dfrac{z}{t},\ \ P-\dfrac{z}{t}
    \end{array}
  $
  
  \item[$\begin{array}{ll}
  4.26\\
  c=0\\ b=1
  \end{array}$]
 
    $\begin{array}{ll}
    \text{Basis}\colon & X_1,\ \ X_4,\ \ aX_5+X_{11},\ \ Y_{1}+X_5,\ \ a\neq0  \\ 
    \text{Inv.:}& \dfrac{z}{t},\ \ v-\dfrac{y}{t},\ \ w,\ \ P+a\ln|t|-\dfrac{y}{t}
    \end{array}
  $
  
   \item[$\begin{array}{ll}
  4.28
  \end{array}$]
 
    $\begin{array}{ll}
    \text{Basis}\colon & X_1,\ \ X_4,\ \ X_{11},\ \ Y_{1}+X_5\\ 
    \text{Inv.:}& \dfrac{z}{t},\ \ v-\dfrac{y}{t},\ \ w,\ \ P-\dfrac{y}{t}
    \end{array}
  $

  \item[$\begin{array}{ll}
  4.29\\
  c\neq 0
  \end{array}$]
 
    $\begin{array}{ll}
    \text{Basis}\colon & X_2,\ \ X_3,\ \ aX_4+bX_5+X_{11},\ \ Y_{1}+cX_4+dX_5+eX_6,\\[2mm] & c^2+d^2+e^2=1,\ \ b\neq0\\ 
    \text{Inv.:}& u-\dfrac{x}{t},\ \ v-\dfrac{dx}{ct}+\left(\dfrac{ad}{c}-b\right)\ln|t|,\ \ w-\dfrac{ex}{ct}+\dfrac{ae}{c}\ln|t|,\\[2mm] & P-\dfrac{x}{ct}+\dfrac{a}{c}\ln|t|
    \end{array}
  $

   \item[$\begin{array}{ll}
  4.29\\
  c=0\\ d\neq 0
  \end{array}$]
 
    $\begin{array}{ll}
    \text{Basis}\colon & X_2,\ \ X_3,\ \ aX_4+bX_5+X_{11},\ \ Y_{1}+dX_5+eX_6,\\[2mm] & d^2+e^2=1,\ \ b\neq0\\ 
    \text{Inv.:}& \dfrac{x}{t}-a\ln|t|,\ \ u-\dfrac{x}{t},\ \ w-\dfrac{e}{d}v+\dfrac{be}{d}\ln|t|,\ \ P-\dfrac{v}{d}+\dfrac{b}{d}\ln|t|
    \end{array}
  $

  \item[$\begin{array}{ll}
  4.29\\
  c=0\\ d=0\\ e=1
  \end{array}$]
 
    $\begin{array}{ll}
    \text{Basis}\colon & X_2,\ \ X_3,\ \ aX_4+bX_5+X_{11},\ \ Y_{1}+X_6,\ \ b\neq0\\ 
    \text{Inv.:}& \dfrac{x}{t}-a\ln|t|,\ \ u-\dfrac{x}{t},\ \ v-b\ln|t|,\ \ P-w
    \end{array}
  $

  \item[$\begin{array}{ll}
  4.31\\
  b\neq 0
  \end{array}$]
 
    $\begin{array}{ll}
    \text{Basis}\colon & X_2,\ \ X_3,\ \ aX_4+X_{11},\ \ Y_{1}+bX_4+cX_5,\ \ b^2+c^2=1\\ 
    \text{Inv.:}&  u-\dfrac{x}{t},\ \ v-\dfrac{cx}{bt}+\dfrac{ac}{b}\ln|t|,\ \ w,\ \ P-\dfrac{x}{bt}+\dfrac{a}{b}\ln|t|
    \end{array}
  $

  \item[$\begin{array}{ll}
  4.31\\
  b=0\\ c=1
  \end{array}$]
 
    $\begin{array}{ll}
    \text{Basis}\colon & X_2,\ \ X_3,\ \ aX_4+X_{11},\ \ Y_{1}+X_5\\ 
    \text{Inv.:}& \dfrac{x}{t}-a\ln|t|,\ \ u-\dfrac{x}{t},\ \ w,\ \ P-v
    \end{array}
  $

  \item[$\begin{array}{ll}
  4.33\\
  b\neq 0
  \end{array}$]
 
    $\begin{array}{ll}
    \text{Basis}\colon & X_1,\ \ X_2+X_4,\ \ X_{10},\ \ Y_{1}+aX_2+bX_3,\ \ a^2+b^2=1\\ 
    \text{Inv.:}& u-y+\dfrac{a}{b}z,\ \ v,\ \ w,\ \ P-\dfrac{z}{b}
    \end{array}
  $
   \item[$\begin{array}{ll}
  4.33\\
  a=1\\ b=0
  \end{array}$]
 
    $\begin{array}{ll}
    \text{Basis}\colon & X_1,\ \ X_2+X_4,\ \ X_{10},\ \ Y_{1}+X_2\\ 
    \text{Inv.:}&  z,\ \ u+P-y,\ \ v,\ \ w
    \end{array}
  $

  \item[$\begin{array}{ll}
  4.36\\
  b\neq 0
  \end{array}$]
 
    $\begin{array}{ll}
    \text{Basis}\colon & X_2,\ \ X_3,\ \ X_4+aX_5+X_{10},\ \ Y_{1}+bX_1+cX_5+dX_6,\\[2mm] & b^2+c^2+d^2=1,\ \ a\neq0\\ 
    \text{Inv.:}& u-t,\ \ v-\dfrac{c}{b}x+\dfrac{c}{2b}t^2-at,\ \ w-\dfrac{d}{b}x+\dfrac{d}{2b}t^2,\ \ P-\dfrac{x}{b}+\dfrac{t^2}{2b}
    \end{array}
  $

  \item[$\begin{array}{ll}
  4.36\\
  b=0\\ c\neq 0
  \end{array}$]
 
    $\begin{array}{ll}
    \text{Basis}\colon & X_2,\ \ X_3,\ \ X_4+aX_5+X_{10},\ \ Y_{1}+cX_5+dX_6,\\[2mm] & c^2+d^2=1, a\neq0\\ 
    \text{Inv.:}& x-\dfrac{t^2}{2},\ \ u-t,\ \ w-\dfrac{d}{c}v+\dfrac{ad}{c}t,\ \ P-\dfrac{v}{c}+\dfrac{a}{c}t
    \end{array}
  $

  \item[$\begin{array}{ll}
  4.36\\
  b=0\\ c=0\\ d=1
  \end{array}$]
 
    $\begin{array}{ll}
    \text{Basis}\colon & X_2,\ \ X_3,\ \ X_4+aX_5+X_{10},\ \ Y_{1}+X_6,\ \ a\neq0\\ 
    \text{Inv.:}& x-\dfrac{t^2}{2},\ \ u-t,\ \ v-at,\ \ P-w
    \end{array}
  $

  \item[$\begin{array}{ll}
  4.37\\
  a\neq 0
  \end{array}$]
 
    $\begin{array}{ll}
    \text{Basis}\colon & X_2,\ \ X_3,\ \ X_5+X_{10},\ \ Y_{1}+aX_1+bX_5+cX_6,\ \ a^2+b^2+c^2=1\\ 
    \text{Inv.:}& u,\ \ v-\dfrac{b}{a}x-t,\ \ w-\dfrac{c}{a}x,\ \ P-\dfrac{x}{a}
    \end{array}
  $
   \item[$\begin{array}{ll}
  4.37\\
  a=0\\ b\neq 0
  \end{array}$]
 
    $\begin{array}{ll}
    \text{Basis}\colon & X_2,\ \ X_3,\ \ X_5+X_{10},\ \ Y_{1}+bX_5+cX_6,\ \ b^2+c^2=1\\ 
    \text{Inv.:}& x,\ \ u,\ \ w+\dfrac{c}{b}(t-v),\ \ P+\dfrac{1}{b}(t-v)
    \end{array}
  $

  \item[$\begin{array}{ll}
  4.37\\
  a=0\\ b=0\\ c=1
  \end{array}$]
 
    $\begin{array}{ll}
    \text{Basis}\colon & X_2,\ \ X_3,\ \ X_5+X_{10},\ \ Y_{1}+X_6\\ 
    \text{Inv.:}& x,\ \ u,\ \ v-t,\ \ P-w
    \end{array}
  $

  \item[$\begin{array}{ll}
  4.39\\
  a\neq 0
  \end{array}$]
 
    $\begin{array}{ll}
    \text{Basis}\colon & X_2,\ \ X_3,\ \ X_4+X_{10},\ \ Y_{1}+aX_1+bX_5,\ \ a^2+b^2=1\\ 
    \text{Inv.:}& u-t,\ \ v-\dfrac{b}{a}x+\dfrac{b}{2a}t^2,\ \ w,\ \ P-\dfrac{x}{a}+\dfrac{t^2}{2a}
    \end{array}
  $

  \item[$\begin{array}{ll}
  4.39\\
  a=0\\ b=1
  \end{array}$]
 
    $\begin{array}{ll}
    \text{Basis}\colon & X_2,\ \ X_3,\ \ X_4+X_{10},\ \ Y_{1}+X_5\\ 
    \text{Inv.:}& x-\dfrac{t^2}{2},\ \ u-t,\ \ w,\ \ P-v
    \end{array}
  $

  \item[$\begin{array}{ll}
  4.40\\
  a\neq 0
  \end{array}$]
 
    $\begin{array}{ll}
    \text{Basis}\colon & X_2,\ \ X_3,\ \ X_{10},\ \ Y_{1}+aX_1+bX_5,\ \ a^2+b^2=1\\ 
    \text{Inv.:}& u,\ \ v-\dfrac{b}{a}x,\ \ w,\ \ P-\dfrac{x}{a}
    \end{array}
  $
   \item[$\begin{array}{ll}
  4.40\\
  a=0\\ b=1
  \end{array}$]
 
    $\begin{array}{ll}
    \text{Basis}\colon & X_2,\ \ X_3,\ \ X_{10},\ \ Y_{1}+X_5\\ 
    \text{Inv.:}& x,\ \ u,\ \ w,\ \ P-v
    \end{array}
  $

  \item[$\begin{array}{ll}
  4.43\\
  c\neq 0
  \end{array}$]
 
    $\begin{array}{ll}
    \text{Basis}\colon & aX_1+X_4,\ \ bX_3+X_5,\ \ bX_2-X_6,\ \ Y_{1}+cX_1+dX_2,\\[2mm] & a^2+b^2=1,\ \ c^2+d^2=1\\[2mm]
    \text{Inv.:}& t,\ \ (t+a)u+cP-x,\\[2mm] & -dt(t+a)u+(t^2+b^2)cv+dtx-cty-bcz,\\[2mm] & bd(t+a)u+(t^2+b^2)cw-bdx+bcy-ctz
    \end{array}
  $

  \item[$\begin{array}{ll}
  4.43\\
  c=0\\ d=1
  \end{array}$]
 
    $\begin{array}{ll}
    \text{Basis}\colon & aX_1+X_4,\ \ bX_3+X_5,\ \ bX_2-X_6,\ \ Y_{1}+X_2,\ \ a^2+b^2=1\\ 
    \text{Inv.:}& t,\ \ u-\dfrac{x}{t+a},\ \ w+\dfrac{b}{t}v-\dfrac{z}{t},\ \ P+\left(\dfrac{b^2}{t}+t\right)v-y-b\dfrac{z}{t}
    \end{array}
  $

  \item[$\begin{array}{ll}
  4.46
  \end{array}$]
 
    $\begin{array}{ll}
    \text{Basis}\colon & X_4,\ \ X_5,\ \ X_6,\ \ Y_{1}+X_1\\ 
    \text{Inv.:}& t,\ \ tu+P-x,\ \ v-\dfrac{y}{t},\ \ w-\dfrac{z}{t}
    \end{array}
  $

  \item[$\begin{array}{ll}
  4.49\\
  a\neq 0
  \end{array}$]
 
    $\begin{array}{ll}
    \text{Basis}\colon & aX_1+X_3,\ \ X_5,\ \ X_6,\ \ Y_{1}+bX_4+X_{11}\\ 
    \text{Inv.:}& u-b\ln|t|,\ \ v-\dfrac{y}{t},\ \ w+\dfrac{x}{at}-\dfrac{z}{t}-\dfrac{b}{a}\ln|t|,\ \ P-\ln|t|
    \end{array}
  $
   \item[$\begin{array}{ll}
  4.49\\
  a=0
  \end{array}$]
 
    $\begin{array}{ll}
    \text{Basis}\colon & X_3,\ \ X_5,\ \ X_6,\ \ Y_{1}+bX_4+X_{11}\\ 
    \text{Inv.:}& \dfrac{x}{t}-b\ln|t|,\ \ u-\dfrac{x}{t},\ \ v-\dfrac{y}{t},\ \ P-\ln|t|
    \end{array}
  $

  \item[$\begin{array}{ll}
  4.52
  \end{array}$]
 
    $\begin{array}{ll}
    \text{Basis}\colon & X_1,\ \ X_3+X_5,\ \ aX_2+X_6,\ \ Y_{1}+bX_2+cX_3+X_4,\ \ a\neq-1\\ 
    \text{Inv.:}& t,\ \ \dfrac{bt-ac}{t^2-a}u+v+\dfrac{az-ty}{t^2-a},\ \ \dfrac{ct-b}{t^2-a}u+w+\dfrac{y-tz}{t^2-a},\ \ P-u
    \end{array}
  $

  \item[$\begin{array}{rr}
  4.53\\
  {\scriptstyle a^2+b^2\neq 0}
  \end{array}$]
 
    $\begin{array}{ll}
    \text{Basis}\colon & X_1,\ \ X_3+X_5,\ \ aX_2+X_6,\ \ Y_{1}+bX_2+cX_3,\ \ b^2+c^2=1\\ 
    \text{Inv.:}& t,\ \ u,\ \ w+\dfrac{1}{ac-bt}((ct-b)v-cy+bz),\\[2mm] & P+\dfrac{1}{ac-bt}((-t^2+a)v+ty-az)
    \end{array}
  $

  \item[$\begin{array}{ll}
  4.53\\
  a=0\\ b=0\\ c=1
  \end{array}$]
 
    $\begin{array}{ll}
    \text{Basis}\colon & X_1,\ \ X_3+X_5,\ \ X_6,\ \ Y_{1}+X_3\\ 
    \text{Inv.:}& t,\ \ u,\ \ v-\dfrac{y}{t},\ \ P+tw+\dfrac{y}{t}-z
    \end{array}
  $

  \item[$\begin{array}{ll}
  4.55
  \end{array}$]
 
    $\begin{array}{ll}
    \text{Basis}\colon & X_1,\ \ X_3+X_5,\ \ X_2-X_6,\ \ Y_{1}+bX_2+cX_3+X_4\\ 
    \text{Inv.:}& t,\ \ v+\dfrac{1}{t^2+1}[(bt+c)u-ty-z],\\[2mm] & w+\dfrac{1}{t^2+1}[(ct-b)u+y-tz],\ \ P-u
    \end{array}
  $
   \item[$\begin{array}{ll}
  4.58
  \end{array}$]
 
    $\begin{array}{ll}
    \text{Basis}\colon & X_1,\ \ X_5,\ \ X_6,\ \ Y_{1}+aX_2+bX_3+X_4,\ \ a^2+b^2=\varepsilon,\ \ \varepsilon=0 \text{ or } 1\\ 
    \text{Inv.:}& t,\ \ v+\dfrac{a}{t}u-\dfrac{y}{t},\ \ w+\dfrac{b}{t}u-\dfrac{z}{t},\ \ P-u
    \end{array}
  $

  \item[$\begin{array}{ll}
  4.59\\
  a\neq 0
  \end{array}$]
 
    $\begin{array}{ll}
    \text{Basis}\colon & X_1,\ \ X_5,\ \ X_6,\ \ Y_{1}+aX_2+bX_3,\ \ a^2+b^2=1\\ 
    \text{Inv.:}& t,\ \ u,\ \ w-\dfrac{b}{a}v+\dfrac{by}{at}-\dfrac{z}{t},\ \ P+\dfrac{t}{a}v-\dfrac{y}{a}
    \end{array}
  $

  \item[$\begin{array}{ll}
  4.59\\
  a=0\\ b=1
  \end{array}$]
 
    $\begin{array}{ll}
    \text{Basis}\colon & X_1,\ \ X_5,\ \ X_6,\ \ Y_{1}+X_3\\ 
    \text{Inv.:}& t,\ \ u,\ \ v-\dfrac{y}{t},\ \ tw+P-z
    \end{array}
  $

  \item[$\begin{array}{ll}
  4.63
  \end{array}$]
 
    $\begin{array}{ll}
    \text{Basis}\colon & X_2,\ \ aX_1+X_3,\ \ X_4,\ \ Y_{1}+X_1\\ 
    \text{Inv.:}& t,\ \ v,\ \ w,\ \ P+tu-x+az
    \end{array}
  $

  \item[$\begin{array}{ll}
  4.66\\
  \varepsilon =1  
  \end{array}$]
 
    $\begin{array}{ll}
    \text{Basis}\colon & X_2,\ \ X_3,\ \ X_4,\ \ Y_{1}+X_1+aX_5+bX_6,\ \ a^2+b^2=1\\[2mm] 
    \text{Inv.:}& t,\ \ tu+P-x,\ \ (tu-x)a+v,\ \ (tu-x)b+w
    \end{array}
  $
   \item[$\begin{array}{ll}
  4.66\\
  a\neq 0\\ \varepsilon=0
  \end{array}$]
 
    $\begin{array}{ll}
    \text{Basis}\colon & X_2,\ \ X_3,\ \ X_4,\ \ Y_{1}+aX_5+bX_6,\ \ a^2+b^2=1\\ 
    \text{Inv.:}& t,\ \ u-\dfrac{x}{t},\ \ w-\dfrac{b}{a}v,\ \ P-\dfrac{v}{a}
    \end{array}
  $

  \item[$\begin{array}{ll}
  4.66\\
  a=0\\ \varepsilon=0\\ b=1
  \end{array}$]
 
    $\begin{array}{ll}
    \text{Basis}\colon & X_2,\ \ X_3,\ \ X_4,\ \ Y_{1}+X_6 \\ 
    \text{Inv.:}& t,\ \ u-\dfrac{x}{t},\ \ v,\ \ P-w
    \end{array}
  $

  \item[$\begin{array}{ll}
  4.67
  \end{array}$]
 
    $\begin{array}{ll}
    \text{Basis}\colon & X_1,\ \ X_2,\ \ X_3+X_4,\ \ Y_{1}+aX_5+bX_6+X_{10}\\ 
    \text{Inv.:}& u-z+\dfrac{b}{2}t^2,\ \ v-at,\ \ w-bt,\ \ P-t
    \end{array}
  $

  \item[$\begin{array}{ll}
  4.68
  \end{array}$]
 
    $\begin{array}{ll}
    \text{Basis}\colon & X_1,\ \ X_2,\ \ X_3+X_4,\ \ Y_{1}+aX_5+X_6\\ 
    \text{Inv.:}& t,\ \ v+a\dfrac{u}{t}-a\dfrac{z}{t},\ \ w+\dfrac{u}{t}-\dfrac{z}{t},\ \ P-w
    \end{array}
  $

  \item[$\begin{array}{ll}
  4.69\\
  a\neq 0
  \end{array}$]
 
    $\begin{array}{ll}
    \text{Basis}\colon & X_1,\ \ X_2,\ \ X_3+X_4,\ \ Y_{1}+aX_3+X_5\\ 
    \text{Inv.:}& t,\ \ v+\dfrac{u}{a}-\dfrac{z}{a},\ \ w,\ \ P-v
    \end{array}
  $
   \item[$\begin{array}{ll}
  4.69\\
  a=0
  \end{array}$]
 
    $\begin{array}{ll}
    \text{Basis}\colon & X_1,\ \ X_2,\ \ X_3+X_4,\ \ Y_{1}+X_5\\ 
    \text{Inv.:}& t,\ \ u-z,\ \ w,\ \ P-v
    \end{array}
  $

  \item[$\begin{array}{ll}
  4.70
  \end{array}$]
 
    $\begin{array}{ll}
    \text{Basis}\colon & X_1,\ \ X_2,\ \ X_3+X_4,\ \ Y_{1}+X_3\\ 
    \text{Inv.:}& t,\ \ v,\ \ w,\ \ P+u-z
    \end{array}
  $

  \item[$\begin{array}{ll}
  4.72
  \end{array}$]
 
    $\begin{array}{ll}
    \text{Basis}\colon & X_1,\ \ X_2,\ \ X_4,\ \ Y_{1}+bX_5+X_6\\ 
    \text{Inv.:}& t,\ \ v-b\dfrac{z}{t},\ \ w-\dfrac{z}{t},\ \ P-\dfrac{z}{t}
    \end{array}
  $

  \item[$\begin{array}{ll}
  4.73\\
  a\neq 0
  \end{array}$]
 
    $\begin{array}{ll}
    \text{Basis}\colon & X_1,\ \ X_2,\ \ X_4,\ \ Y_{1}+aX_3+bX_5,\ \ a^2+b^2=1\\ 
    \text{Inv.:}& t,\ \ v-\dfrac{b}{a}z,\ \ w,\ \ P-\dfrac{z}{a}
    \end{array}
  $

  \item[$\begin{array}{ll}
  4.73\\
  a=0\\ b=1
  \end{array}$]
 
    $\begin{array}{ll}
    \text{Basis}\colon & X_1,\ \ X_2,\ \ X_4,\ \ Y_{1}+X_5\\ 
    \text{Inv.:}& t,\ \ z,\ \ w,\ \ P-v
    \end{array}
  $
   \item[$\begin{array}{ll}
  4.75
  \end{array}$]
 
    $\begin{array}{ll}
    \text{Basis}\colon & X_1,\ \ X_2,\ \ X_3,\ \ Y_{1}+aX_4+X_{11}\\ 
    \text{Inv.:}& u-a\ln|t|,\ \ v,\ \ w,\ \ P-\ln|t|
    \end{array}
  $

  \item[$\begin{array}{ll}
  4.76
  \end{array}$]
 
    $\begin{array}{ll}
    \text{Basis}\colon & X_1,\ \ X_2,\ \ X_3,\ \ Y_{1}+\varepsilon X_4+X_{10}\\ 
    \text{Inv.:}& u-\varepsilon t,\ \ v,\ \ w,\ \ P-t
    \end{array}
  $

    \end{enumerate}
\setlength\extrarowheight{8pt}

\section{Submodel}
With a goal of obtaining new exact solutions of ~\eqref{siraeva:eq0}, \eqref{Siraeva:US_special}, I have chosen a subalgebra 4.73, $a\neq 0$, from the previous section because of its simplicity. The basis generators of  subalgebra $4.73$~\cite{siraeva2014optimal} are
\begin{equation}\label{siraeva:eq1}
\begin{array}{c}
X_1=\partial_{x},\ \ X_2=\partial_{y},\ \ X_4=t\partial_{x}+\partial_{u}, \\[2mm]
Y_1+aX_3+bX_5=\partial_P+a\partial_x+bt\partial_{y}+b\partial_{v}, \  \ a^2+b^2=1. 
\end{array}
\end{equation}
Using invariants of subalgebra 4.73, $a\neq 0$, from Section~\ref{siraeva:section3}, we introduce a new set of dependent variables 
\begin{equation}\label{siraeva:eq2}
	\begin{array}{c}
u=u(t,x,y,z),\ \ v=\dfrac{b}{a}z+v_1(t),\ \ w=w_1(t),\ \ \rho=\rho(t),\\ P=P_1(t)+\gamma\dfrac{z}{a}, \ \ S=S_1(t)+\gamma\dfrac{z}{a},
	\end{array}
\end{equation}
where $\gamma=1$ in the case of Lie algebra $L_{12}$ and $\gamma=0$ in the case of Lie algebra $L_{11}$. Thus, parameter $\gamma$ allows distinguish the results for 11-- and 12--dimensional Lie algebras.

Substitution \eqref{siraeva:eq2} in \eqref{siraeva:eq0}, \eqref{Siraeva:US_special} leads to a partially invariant submodel of rank 1 and defect 1
\begin{equation}\label{siraeva:eq3}
	\begin{array}{c}
u_t+uu_x+\left(\dfrac{b}{a}z+v_1\right)u_y+w_1u_z=0,\\[2mm]
v_{1t}+\dfrac{b}{a}w_1=0,\\[2mm]
w_{1t}+\dfrac{\gamma}{a}\rho^{-1}=0,\\[2mm]
\rho_t+\rho u_x=0,\\[2mm]
S_{1t}+\dfrac{\gamma}{a}w_1=0,\ \ P_1=S_1+f(\rho).
	\end{array}
\end{equation}
\subsection{Solution for isochoric motion of media}
If $\rho=\rho_0=\mathrm{const}$, exact particular solution of system \eqref{siraeva:eq0}, \eqref{Siraeva:US_special} from submodel \eqref{siraeva:eq3} has the following form due to transformations $4^o$ from $\eqref{Sir:GroupGal_rashir}$ with $\vec{b}=(0,-v_0,-w_0)$ and $6^o$ from $\eqref{Sir:peren_po_p}$ with $P_0=-f(\rho_0)-S_0$ 
\begin{equation}\label{siraeva:eq9}
	\begin{array}{c}
		u=\Phi\left(z+\dfrac{\gamma}{2a\rho_0}t^2, y-\dfrac{b}{a}tz-\dfrac{\gamma t^3b}{2a^2\rho_0}\right),\\[2mm]
		v=\dfrac{b}{a}z+\dfrac{\gamma b}{2a^2\rho_0}t^2,\\[2mm]
		w=-\dfrac{\gamma t}{a\rho_0},\\[2mm]
		\rho=\rho_0,\\[2mm]
		P=\dfrac{\gamma}{a}z+\dfrac{\gamma^2}{2a^2\rho_0}t^2,\\[2mm]
		S=\dfrac{\gamma}{a}z+\dfrac{\gamma^2}{2a^2\rho_0}t^2.
	\end{array}
\end{equation}
Particle motion is given by the equation~\cite{Sir:lit:Ovs3}:
\begin{equation}\label{Sir:eq_mirlin}
	\begin{array}{c}
		\dfrac{d\vec{x}}{dt}=\vec{u}(t,\vec{x}).\nonumber
	\end{array}
\end{equation}
The world lines of the particles in $\mathbb{R}^4(t,\vec{x})$ from \eqref{siraeva:eq9} are 
\begin{equation}\label{siraeva:eq20} 
	\begin{array}{c}
x=t\Phi(z_0,y_0)+x_0,\\
y=\dfrac{b}{a}z_0t+y_0,\\
z=-\dfrac{\gamma t^2}{2a\rho_0}+z_0.
	\end{array}
\end{equation}
which are the coordinates of particles at $t=0$
\begin{equation}
	\begin{array}{c}
		x(0)=x_0,\ \ y(0)=y_0,\ \ z(0)=z_0.
	\end{array}
\end{equation}
The Jacobian determinant of the transformation~\eqref{siraeva:eq20} from the Lagrangian to Euler coordinates is
\begin{equation}
	\begin{array}{c}
		J=\begin{vmatrix}
			1& t\Phi_{y_0} & t\Phi_{z_0}\\
			0& 1& \frac{b}{a}t\\
			0& 0& 1
		\end{vmatrix}
	\end{array}=1.
\end{equation}
The world lines of the particles do not intersect. 
Solution \eqref{siraeva:eq9} has no collapse, the motion of media is isochoric. Projections of the world lines \eqref{siraeva:eq20} to $\mathbb{R}^3(\vec{x})$ are particle trajectories.

The motion of particles is vortex if $b\neq 0$ or $\Phi\neq const$
\begin{equation}\label{siraeva:eq21}
	\begin{array}{c}
		\vec{\omega}=\text{rot}\vec{u}=(w_y-v_z,u_z-w_x,v_x-u_y)=\\
  =\left(-b/a,\Phi_{z_0}-\dfrac{b}{a}t\Phi_{y_0},-\Phi_{y_0}\right).
	\end{array}
\end{equation}
Let us rewrite the gas-dynamic functions \eqref{siraeva:eq9} in terms of the Lagrangian coordinates
\begin{equation}\label{siraeva:eq22}
	\begin{array}{c}
	u=\Phi(z_0,y_0), \ \ v=\dfrac{b}{a}z_0,\ \ w=-\gamma\dfrac{t}{a\rho_0} \\
	\rho=\rho_0, \ \ P=\dfrac{\gamma}{a}z_0.
	\end{array}
\end{equation}
The velocity of the particle changes only in the direction of the OZ axis.

Projections of trajectories~\eqref{siraeva:eq20} on plane $(x,y)$ are straight lines with direction vector $\vec{q}=\left(\Phi(z_0,y_0),\dfrac{b}{a}z_0\right)$ with length
$q=|\Vec{q}|=\sqrt{\Phi^2+\dfrac{b^2}{a^2}z_0^2}$.
These straight lines are given by the following formula for $\Phi\neq 0$
\begin{equation}\label{siraeva:eq31}
	\begin{array}{c}
	y=\dfrac{bz_0}{a\Phi}x+y_0-\dfrac{bz_0x_0}{a\Phi}.
	\end{array}
\end{equation}
Let $s=tq$ be the length of the rectilinear path of the particle projection in the plane $(x,y)$. In the plane $(s,z)$ the trajectory of a particle is a parabola (except in the case of $bz_0=\Phi=0$)
\begin{equation}\label{siraeva:eq32}
	\begin{array}{c}
	z=-\dfrac{\gamma}{2a\rho_0}\left(\dfrac{s}{q}\right)^2+z_0.
	\end{array}
\end{equation}
For $\Phi=0$ $(bz_0\neq 0)$ parabolas~\eqref{siraeva:eq20} are in the planes $x=x_0$.

For $\Phi=bz_0=0$ particle trajectories~\eqref{siraeva:eq20} are rays.

\begin{figure}[h!]
        \centering
	\includegraphics[scale=0.5]{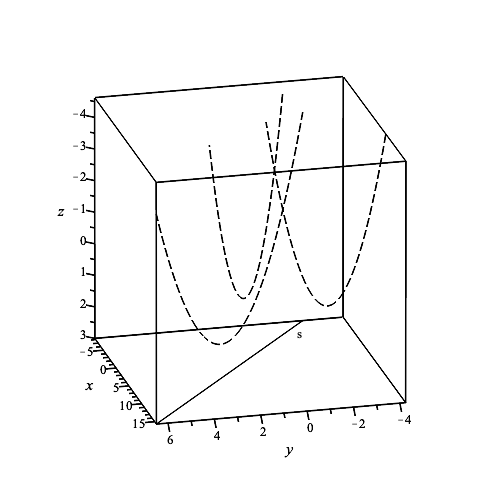}
	\caption{Trajectories of particles~\eqref{siraeva:eq20} with coordinates $(x_0,y_0,z_0)=(5,-2,1)$ and with coordinates $(x_0,y_0,z_0)=(1,1,1);(11/3,2,2)$ according to formulas~~\eqref{siraeva:eq33} at $a=\dfrac{4}{5}$, $b=\dfrac{3}{5}$, $\rho_0=1$, $\Phi=y_0+z_0$, $t=-3..3$.}\label{siraeva:f1}
\end{figure}

For $\Phi\neq 0$ two parabolas~\eqref{siraeva:eq20} with $(x_0,y_0,z_0)=(x^1_0,y^1_0,z^1_0);(x^2_0,y^2_0,z^2_0)$ are in the same plane if
\begin{equation}\label{siraeva:eq33}
	\begin{array}{c}
	\dfrac{z^1_0}{\Phi(z^1_0,y^1_0)}=\dfrac{z^2_0}{\Phi(z^2_0,y^2_0)},\\[4mm]
 y^1_0-\dfrac{bz^1_0x^1_0}{a\Phi(z^1_0,y^1_0)}= y^2_0-\dfrac{bz^2_0x^2_0}{a\Phi(z^2_0,y^2_0)}.
	\end{array}
\end{equation}
Examples of trajectories~\eqref{siraeva:eq20} are shown in Fig.~\ref{siraeva:f1}.

Thus, particles move along parabolas or rays, moving first in one direction of the OZ axis for $t<0$ (depending on a), and then in the other direction of the OZ axis for $t>0$. The surfaces of the level of equal pressure $P=\text{const}$ at each fixed moment of time are planes $z=\text{const}$.

The second and the third equalities of~\eqref{siraeva:eq20} for $b\neq 0$ give a relation
\begin{equation}\label{siraeva:eq25}
	\begin{array}{c}
	y-\dfrac{b}{a}\left(z+\dfrac{\gamma t^2}{2a\rho_0}\right)t=y_0,
	\end{array}
\end{equation}
which defines the moving plane perpendicular to the vector for fixed $t$
\begin{equation}\label{siraeva:eq26}
	\begin{array}{c}
	\Vec{n}_1=\left(0,1,-\dfrac{b}{a}t\right).
	\end{array}
\end{equation}
Therefore, the particles that are on the plane $y=y_0$ at the initial time $t=0$, will be on the plane~\eqref{siraeva:eq25} at any other time $t$.

\begin{figure}[ht]
        \centering
	\includegraphics[scale=0.4]{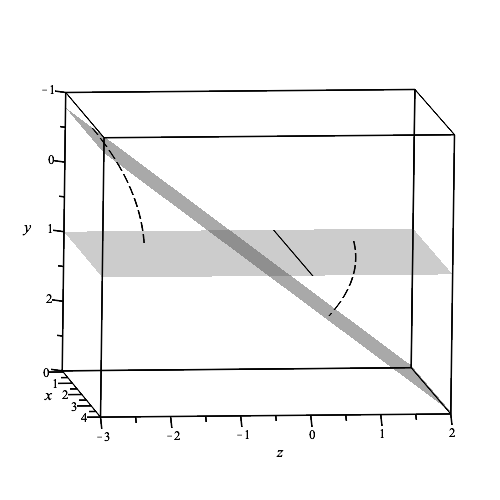}
	\caption{Plane~\eqref{siraeva:eq25} with $y_0=1$ at $t=0$ (light grey color) and $t=1$ (grey color) and trajectories of two particles~\eqref{siraeva:eq20} (dash lines) at $t\in [0,1]$ with coordinates $(x_0, y_0, z_0)=(1, 1, 1);$ $(1;1;-2)$, $\Phi=y_0+z_0$, $a=\dfrac{4}{5}$, $b=\dfrac{3}{5}$, $\rho_0=1$. The solid line is the straight line of intersection of planes~\eqref{siraeva:eq25} and~\eqref{siraeva:eq27}, separating the particles flying out on different sides from the plane~\eqref{siraeva:eq25} at $t=0$.}\label{siraeva:f2}
\end{figure}

The condition $\Vec{u}\cdot\Vec{n}_1=0$ defines the moving plane 
\begin{equation}\label{siraeva:eq27}
	\begin{array}{c}
	z=-\dfrac{3\gamma}{2a\rho_0}t^2,
	\end{array}
\end{equation}
which separates the particles flying out on opposite sides of these planes (Fig~\ref{siraeva:f2}).
\subsection{Solution for non-isochoric motion of media}
If $\rho\neq\mathrm{const}$, exact particular solution of system \eqref{siraeva:eq0}, \eqref{Siraeva:US_special} from submodel \eqref{siraeva:eq3} has the following form due to transformations $4^o$ from $\eqref{Sir:GroupGal_rashir}$ with $\vec{b}=(0,-v_0,-w_0)$ and $6^o$ from $\eqref{Sir:peren_po_p}$ with $P_0=-S_0$ 
\begin{equation}\label{siraeva:eq5}
	\begin{array}{c}
		u=\dfrac{x+\Phi\left(z+\dfrac{\gamma t^3}{6a\rho_0}, y-\dfrac{b}{a}zt-\dfrac{\gamma b}{6a^2\rho_0}t^4\right)}{t},\\[2mm]
		v=\dfrac{b}{a}z+\dfrac{\gamma b}{6a^2\rho_0}t^3,\\[2mm]
		w=-\dfrac{\gamma t^2}{2a\rho_0},\\[2mm]
		\rho=\dfrac{\rho_0}{t},\\[2mm]
		P=\dfrac{\gamma}{a}z+\dfrac{\gamma^2}{6a^2\rho_0}t^3+f\left(\dfrac{\rho_0}{t}\right),\\[2mm]
		S=\dfrac{\gamma}{a}z+\dfrac{\gamma^2}{6a^2\rho_0}t^3.
	\end{array}
\end{equation}
The world lines of particles are
\begin{equation}\label{siraeva:eq6}
	\begin{array}{c}
		x=u_0t-\Phi(z_0, y_0),\\
		y=\dfrac{b}{a}z_0t+y_0,\\
		z=-\dfrac{\gamma t^3}{6a\rho _0}+z_0.
	\end{array}
\end{equation}
where $u_0$, $y_0$, $z_0$ are local Lagrangian coordinates ($u_0$ determines the particle velocity in the direction of the Ox axis, $y_0$ and $z_0$ are the y and z coordinates of the particle at $t=0$).

The Jacobian of the transformation~\eqref{siraeva:eq6} from the Lagrangian to Euler coordinates is
\begin{equation}
	\begin{array}{c}
		J=\begin{vmatrix}
			t& -\Phi_{y_0} &-\Phi_{z_0}\\
			0& 1& \frac{b}{a}t\\
			0& 0& 1
		\end{vmatrix}
	\end{array}=t.
\end{equation}
The rank of the Jacobian matrix of transformations \eqref{siraeva:eq6} at the moment of blow-up $t=0$ is equal to 2. At time $t=0$, the particle is on a second-order surface defined by an arbitrary function $\Phi$
\begin{equation}\label{siraeva:eq8}
	\begin{array}{c}
		x=-\Phi(z_0, y_0),\ \ y=y_0, \ \ z=z_0.
	\end{array}
\end{equation}
At each point $(-\Phi(z_0, y_0), y_0, z_0)$ there is an infinite set of particles that differ from each other in speed $u_0$.

The motion of particles is vortex (except for the case $b=0$, $\Phi=\Phi_0$)
\begin{equation}\label{siraeva:eq23}
	\begin{array}{c}
		\vec{\omega}=\text{rot}\vec{u}=\left(-\dfrac{b}{a},\dfrac{1}{t}\Phi_{z_0}-\dfrac{b}{a}\Phi_{y_0},-\dfrac{1}{t}\Phi_{y_0}\right).
	\end{array}
\end{equation}
Let us rewrite the gas-dynamic functions \eqref{siraeva:eq5} in terms of the Lagrangian coordinates
\begin{equation}\label{siraeva:eq24}
	\begin{array}{c}
	u=u_0, \ \ v=\dfrac{b}{a}z_0,\ \ w=-\dfrac{\gamma}{2a\rho_0}t^2, \\
	\rho=\dfrac{\rho_0}{t}, \ \ P=\dfrac{\gamma}{a}z_0+f\left(\dfrac{\rho_0}{t}\right),\ \ S=\dfrac{\gamma}{a}z_0.
	\end{array}
\end{equation}
The velocity of the particle changes in the direction of the OZ axis. 

Particle trajectories~\eqref{siraeva:eq6} for $u_0=bz_0=0$ are straight lines.

For $u_0^2+(bz_0)^2\neq 0$, the projections of trajectories onto the plane $(x,y)$, similarly to the previous case, are straight lines with a direction vector $\vec{q}=\left(u_0,\dfrac{b}{a}z_0\right)$ of length $q=|\vec{q}|=\sqrt{u^2_0+\dfrac{b^2}{a^2}z^2_0}$.
These lines are given by the equation for $u_0\neq 0$
\begin{equation}
	\begin{array}{c}
	y=\dfrac{bz_0}{au_0}x+y_0+\dfrac{bz_0}{au_0}\Phi.
	\end{array}
\end{equation}
Let $s=tq$ be the length of the rectilinear path of the particle projection in the plane $(x,y)$. In the plane $(s,z)$ the trajectory of a particle is a cubic parabola (except for the case of $u_0=bz_0=0$)
\begin{equation}\label{siraeva:eq}
	\begin{array}{c}
	z=-\dfrac{\gamma}{6a\rho_0}\left(\dfrac{s}{q}\right)^3+z_0.
	\end{array}
\end{equation}
For $u_0=0$ $(bz_0\neq 0)$ parabolas~\eqref{siraeva:eq6} are in planes $x=-\Phi (z_0,y_0)=\text{const}$.

\begin{figure}[H]
        \centering
	\includegraphics[scale=0.35]{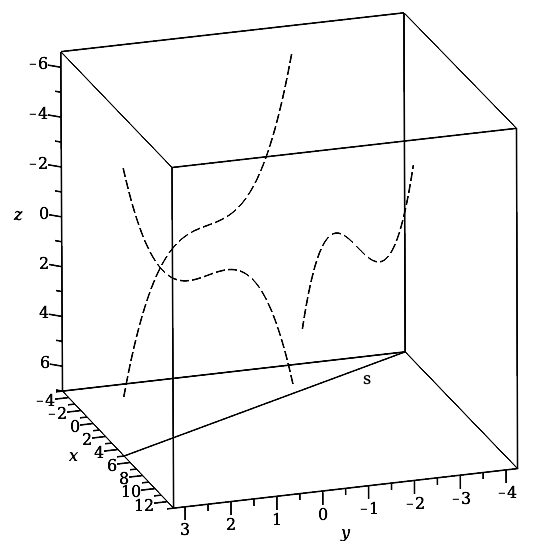}
	\caption{Trajectories of particles~\eqref{siraeva:eq6} with coordinates $(u_0,y_0,z_0)=(3,-2,0)$ and with coordinates $(u_0,y_0,z_0)=(1,1,1);(-1,1,-1)$ according to formulas~~\eqref{siraeva:eq34} at $a=\dfrac{4}{5}$, $b=\dfrac{3}{5}$, $\rho_0=1$, $\Phi=-y_0^2-z_0^2$, $t=-3..3$.}\label{siraeva:f3}
\end{figure}

For $u_0\neq 0$ two parabolas~\eqref{siraeva:eq6} with $(u_0,y_0,z_0)=(u^1_0,y^1_0,z^1_0);(u^2_0,y^2_0,z^2_0)$ are in the same plane if
\begin{equation}\label{siraeva:eq34}
	\begin{array}{c}
	\dfrac{z^1_0}{u_0^1}=\dfrac{z^2_0}{u_0^2},\\[4mm]
 y^1_0+\dfrac{bz^1_0}{au^1_0}\Phi(z^1_0,y^1_0)= y^2_0+\dfrac{bz^2_0}{au_0^2}\Phi(z^2_0,y^2_0).
	\end{array}
\end{equation}
Examples of trajectories~\eqref{siraeva:eq6} are shown in Fig.~\ref{siraeva:f3}.

\begin{figure}[H]
        \centering
	\includegraphics[scale=0.35]{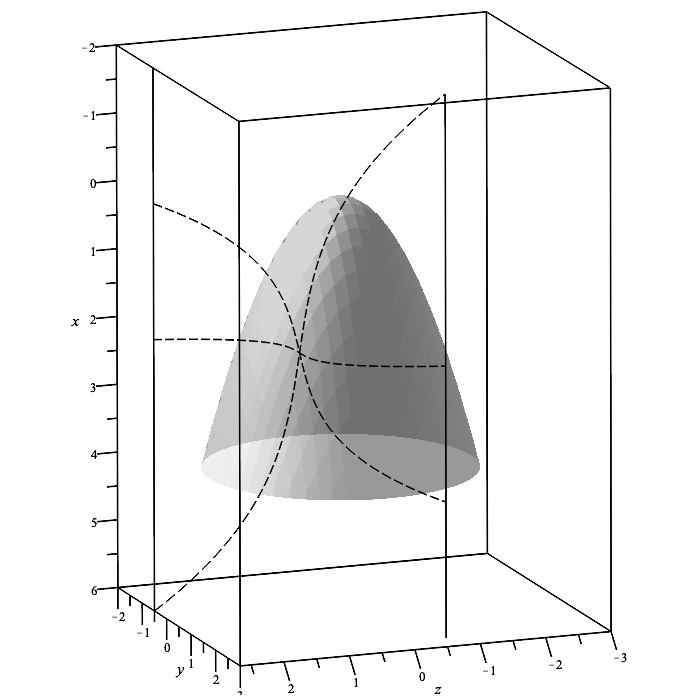}
	\caption{The surface of blow-up~\eqref{siraeva:eq8} in the form of an elliptical paraboloid $\Phi=-y_0^2-z_0^2$, trajectories~\eqref{siraeva:eq6} of three particles with coordinates $(u_0,y_0,z_0)=(1,1,1);(-2,1,1);(0,1,1)$ flying out (flying in) from one point of blow-up and the straight lines on which they end up at the moments of time $t=-2$ and $t=2$; $a=\dfrac{4} {5}$, $b=\dfrac{3}{5}$, $\rho_0=1$, $\Phi=-y_0^2-z_0^2$, $t=-2..2$.}\label{siraeva:f4}
\end{figure}

From one point of blow-up $(-\Phi(z_0,y_0),y_0,z_0)$, particles fly out at different speeds and end up on a straight line at any other fixed moment in time $t\neq 0$. The equation of this line is given by formulas~\eqref{siraeva:eq6} if we consider $t$, $y_0$, $z_0$ fixed in them, and $u_0$ is a parameter. This line is parallel to the Ox axis (Fig~\ref{siraeva:f4}).
It is worth noting that projections of the trajectories onto the plane $(y,z)$ of all points flying from one point of blow-up coincide and represent a cubic parabola given by the second and third equations of formulas~\eqref{siraeva:eq6}.
\section{Conclusion}
In this article, four--dimensional subalgebras of a 12--dimensional Lie algebra, admitted by the equations of gas dynamics with a special equation of state, have been considered for the first time. We have calculated invariants for 50 four--dimensional subalgebras. We have used invariants from one subalgebra and  have calculated a partially invariant submodel of rank one and defect one, a system of partial differential equations. The submodel leads to two families of new exact solutions of the gas dynamics equations, one of which describes the isochoric motion of media, and the other solution specifies an arbitrary blow-up surface. For the first family of solutions, the particle trajectories are parabolas or rays; for the second family of solutions, the particles move along cubic parabolas or straight lines. From each point of the blow-up surface, particles fly out at different speeds and end up on a straight line at any other fixed moment of time. A description of the motion of particles for each family of solutions have been given.

\bibliographystyle{unsrt}  
\bibliography{references}  

\begin{thebibliography}{10}

\bibitem{Ovs1}
L.V. Ovsyannikov.
\newblock The <<podmodeli>> program. gas dynamics.
\newblock {\em Journal of Applied Mathematics and Mechanics}, 1994.
\newblock
  \href{https://doi.org/10.1016/0021-8928(94)90137-6}{doi.org/10.1016/0021-8928(94)90137-6}.

\bibitem{Sir:lit:Ovs2}
L.V. Ovsiannikov.
\newblock {\em Group analysis of differential equations}.
\newblock Academic press, 1982.
\newblock
  \href{https://www.sciencedirect.com/book/9780125316804/group-analysis-of-differential-equations}{www.sciencedirect.com}.

\bibitem{Sir:Olver1986}
P.J. Olver.
\newblock {\em Applications of Lie groups to differential equations}.
\newblock New York: Springer, 1986.
\newblock
  \href{https://doi.org/10.1007/978-1-4684-0274-2}{doi.org/10.1007/978-1-4684-0274-2}.

\bibitem{Sir:Bluman1989}
G.W. Bluman and S.~Kumei.
\newblock {\em Symmetries and differential equations}.
\newblock New York: Springer, 1989.
\newblock
  \href{https://doi.org/10.1007/978-1-4757-4307-4}{doi.org/10.1007/978-1-4757-4307-4}.

\bibitem{Sir:Ibrag1994}
N.Kh Ibragimov.
\newblock {\em CRC Handbook of Lie Group Analysis of Differential Equations,
  Volume I: Symmetries, Exact Solutions, and Conservation Laws}.
\newblock CRC Press, Boca Raton, 1994.

\bibitem{Sir:Ibrag2001}
N.H. Ibragimov.
\newblock {\em Transformation groups applied to mathematical physics}.
\newblock Springer Netherlands, 2001.

\bibitem{Sir:Hydon2000}
P.E. Hydon.
\newblock {\em Symmetry Methods for Differential Equations: A Beginner's
  Guide}.
\newblock New York: Cambridge University Press, 2000.

\bibitem{Sir:Arrigo2015}
D.J. Arrigo.
\newblock {\em Symmetry Analysis of Differential Equations: An Introduction}.
\newblock USA: Wiley, 2015.

\bibitem{siraeva2014optimal}
D.T. Siraeva.
\newblock Optimal system of non-similar subalgebras of sum of two ideals.
\newblock {\em Ufa Mathematical Journal}, 6(1):90--103, 2014.
\newblock
  \href{https://doi.org/10.13108/2014-6-1-90}{doi.org/10.13108/2014-6-1-90}.

\bibitem{SirMnSys}
D.T. Siraeva.
\newblock Reduction of partially invariant submodels of rank 3 defect 1 to
  invariant submodels.
\newblock {\em Multiphase Systems}, 13(3):59--63, 2018.
\newblock
  \href{https://doi.org/10.21662/mfs2018.3.009}{doi.org/10.21662/mfs2018.3.009}.

\bibitem{Sir16}
D.T. Siraeva.
\newblock Classification of rank 2 stationary submodels of ideal hydrodynamics.
\newblock {\em Chelyabinsk Physical and Mathematical Journal}, 4(1):18--32,
  2019.
\newblock
  \href{https://doi.org/10.24411/2500-0101-2019-14102}{doi.org/10.24411/2500-0101-2019-14102}.

\bibitem{SirSibJIM1}
D.T. Siraeva.
\newblock The canonical form of the rank 2 invariant submodels of evolutionary
  type in ideal hydrodynamics.
\newblock {\em Journal of Applied and Industrial Mathematics}, 13(2):340--349,
  2019.
\newblock
  \href{https://doi.org/10.1134/S1990478919020157}{doi.org/10.1134/S1990478919020157}.

\bibitem{Sir4Khab5}
D.T. Siraeva.
\newblock Invariant submodel of rank 2 on subalgebra of translations linear
  combinations for a hydrodynamic type model.
\newblock {\em Chelyabinsk Physical and Mathematical Journal}, 3(1):38--57,
  2018.
\newblock (in Russian)
  \href{https://cpmj.csu.ru/index.php/cpmj/article/view/133/114}{cpmj.csu.ru}.

\bibitem{Sir:lit:Sir1JoPh20}
D.T. Siraeva.
\newblock Two invariant submodels of rank 1 of the hydrodynamic type equations
  and exact solutions.
\newblock {\em Journal of Physics: Conference Series}, 1666(012049), 2020.
\newblock
  \href{https://doi.org/10.1088/1742-6596/1666/1/012049}{doi.org/10.1088/1742-6596/1666/1/012049}.

\bibitem{Sir:lit:SEMR2021}
D.T. Siraeva.
\newblock Invariant solutions of the gas dynamics equations from 4-parameter
  three-dimensional subalgebras containing all translations in space and
  pressure translation.
\newblock {\em Siberian Electronic Mathematical Reports}, 18(2):1639--1650,
  2021.
\newblock
  \href{https://doi.org/10.33048/semi.2021.18.123}{doi.org/10.33048/semi.2021.18.123}.

\bibitem{Sir:lit:MDPI2022}
R.~Nikonorova, D.~Siraeva, and Y.~Yulmukhametova.
\newblock New exact solutions with a linear velocity field for the gas dynamics
  equations for two types of state equations.
\newblock {\em Mathematics}, 10(1), 2022.
\newblock
  \href{https://doi.org/10.3390/math10010123}{doi.org/10.3390/math10010123}.

\bibitem{Sir:lit:Khabirov2019}
S.V. Khabirov.
\newblock Simple partially invariant solutions.
\newblock {\em Ufa Mathematical Journal}, 11(1):90--99, 2019.
\newblock
  \href{https://doi.org/10.13108/2019-11-1-90}{doi.org/10.13108/2019-11-1-90}.

\bibitem{Sir:lit:Nikon2023}
R.F. Nikonorova.
\newblock Simple invariant solutions of the dynamic equation for a monatomic
  gas.
\newblock {\em Proceedings of the Steklov Institute of Mathematics},
  321(1):S186–S203, 2023.
\newblock
  \href{https://doi.org/10.1134/S0081543823030161}{doi.org/10.1134/S0081543823030161}.

\bibitem{Khab1}
S.V. Khabirov.
\newblock Nonisomorphic lie algebras admitted by gasdynamic models.
\newblock {\em Ufa Mathematical Journal}, 3(2):85--88, 2011.
\newblock \href{https://matem.anrb.ru/en/article?art_id=67}{matem.anrb.ru}.

\bibitem{Khab2}
S.V. Khabirov.
\newblock Optimal system for sum of two ideals admitted by hydrodynamic type
  equations.
\newblock {\em Ufa Mathematical Journal}, 6(2):97--101, 2014.
\newblock
  \href{https://doi.org/10.13108/2014-6-2-97}{doi.org/10.13108/2014-6-2-97}.

\bibitem{Chir1Khab4}
Yu.A. Chirkunov and S.V. Khabirov.
\newblock {\em Elements of Symmetry Analysis of Differential Equations of
  Continuum Mechanics: monograph}.
\newblock NSTU, 2012.

\bibitem{Sir:lit:Khab3}
S.V. Khabirov.
\newblock {\em Lectures analytical methods in gas dynamics}.
\newblock BSU, 2013.

\bibitem{Sir:lit:Mamont1}
E.V. Mamontov.
\newblock Invariant submodels of rank two of the equations of gas dynamics.
\newblock {\em Journal of Applied Mechanics and Technical Physics},
  40(2):232--237, 1999.
\newblock
  \href{https://doi.org/10.1007/BF02468519}{doi.org/10.1007/BF02468519}.

\bibitem{Sir:lit:Mamont2}
E.V. Mamontov.
\newblock Gruppovyye svoystva 2-podmodeley klassa s uravneniy gazovoy dinamiki
  [group properties of 2-submodels of class s of equations of gas dynamics].
\newblock {\em Vestnik Novosibirskogo gosudarstvennogo universiteta. Seriya:
  Matematika, mekhanika, informatika}, 7(1):72--84, 2007.

\bibitem{Sir:lit:Ovs3}
L.V. Ovsiannikov.
\newblock {\em Lectures on the foundations of gas dynamics}.
\newblock 2 Eds., Moscow-Izhevsk: Institut komp'yuternykh issledovaniy, 2003.

\end{thebibliography}

\end{document}